\documentclass[a4paper,10pt,
colorlinks,urlcolor=black,linkcolor=black,citecolor=black,
]{article}
\usepackage{amsmath,amssymb}
\usepackage[dvips]{color}
\usepackage{ascmac}
\usepackage{cite}
\usepackage{hangcaption}
\usepackage{amsmath,amssymb}
\usepackage[dvips]{color}
\usepackage{cite}
\usepackage{eepic}
\usepackage{epic}
\usepackage{hangcaption}
\usepackage{amsmath}
\usepackage{amssymb}
\usepackage{amscd}
\usepackage{amsthm}
\usepackage{multicol}
\usepackage[dvipdfmx,%
 bookmarks=true,%
 bookmarksnumbered=true,%
 setpagesize=false,%
 pdfkeywords={TeX; dvipdfmx; hyperref; color;}]{hyperref}
\usepackage{hyperref}
\def\cases{\left\{\begin{array}{ll}}
\def\endcases{\end{array}\right.}

\def\bigtimes{\mathop{\mbox{\Large $\times$}}}
\begin{document}
\setcounter{page}{1}
\vskip1.5cm
\begin{center}
{\LARGE \bf 
Regression analysis in quantum language
}
\vskip0.5cm
{\bf
\large
Shiro Ishikawa
}
\\
\vskip0.2cm
\rm
\it
Department of Mathematics, Faculty of Science and Technology,
Keio University, 
\\
3-14-1,Hiyoshi, Kohoku-ku, Yokohama,
Japan.
E-mail:
ishikawa@math.keio.ac.jp
\end{center}
\par
\rm
\vskip0.3cm
\par
\noindent
{\bf Abstract}
\normalsize
\vskip0.5cm
\par
\noindent
Although regression analysis has a great history, we consider that
it has always continued being confused.
For example,
the fundamental terms in regression analysis
(e.g.,
"regression", "least-squares method", "explanatory variable", "response variable", etc.)
seem to be historically conventional, that is, these words do not express the essence of regression analysis.
Recently, we proposed quantum language (or, classical and quantum measurement theory), which is characterized as the linguistic turn of the Copenhagen interpretation of quantum mechanics. We believe that this language has a great power of
description,
and therefore, even statistics can be described by quantum language.
Therefore,
in this paper,
we discuss the regression analysis and the generalized linear model
(i.e.,
multiple regression analysis)
in quantum language,
and clarify that the terms
"explanatory variable"
and
"response variable"
is respectively characterized as a kind of causality and the measured value.
%
\par
\vskip1.0cm
\par

\par

\def\Cal{\cal}
\def\bigstimes{\text{\large $\: \boxtimes \,$}}

\par
\noindent

\vskip0.2cm
\par
\noindent
\par
\noindent
\section{
Introduction
}
%

\rm
\par
\par
\noindent
\subsection{The least-squared method in applied mathematics}
\label{subsec:11Least}
\par
\noindent
\par
Let us start from the simple explanation of the least-squared method.
Let
$\{ ({a}_{i}, x_{i} )\}_{{i}=1}^{n}$ be a sequence in the two dimensional real space
${\mathbb R}^2$. Let
$\phi^{(\beta_1, \beta_2)}: {\mathbb R} \to {\mathbb R}$
be the simple function
such that
${\mathbb R} \ni a \mapsto
x= \phi^{(\beta_1, \beta_2)}({a})$
$=\beta_1 {a} + \beta_0 \in {\mathbb R}$,
where
the pair
$(\beta_1, \beta_2) (\in {\mathbb R}^2 )$
is assumed to be unknown.
Define the error $\sigma_{}^{}$
by
\begin{align}
\sigma_{}^2
(\beta_1, \beta_2)=
\frac{1}{{n}} \sum_{{i}=1}^n (x_{i} - \phi^{(\beta_1, \beta_2)}({a}_{i}))^2
\Big(
=
\frac{1}{{n}}
\sum_{{i}=1}^{n}
(
x_{i} -( \beta_1 {a}_{i} + \beta_0 ))^2
\Big)
\end{align}
Then, we have the following minimization problem:
\par
\noindent
\bf
Problem 1
\rm
[The least-squared method].
\begin{itemize}
\item[(A)]
Find the $(\hat{\beta}_0, \hat{\beta}_1)$ 
$( \in {\mathbb R}^2 )$ such that
\begin{align}
\sigma_{}^2
(\hat{\beta}_0, \hat{\beta}_1)
=
\min_{(\beta_1, \beta_2)
\in {\mathbb R}^2
}
\sigma_{}^2
(\beta_1, \beta_2)
\Big(
=
\frac{1}{{n}}
\min_{(\beta_1, \beta_2)
\in {\mathbb R}^2
}
\sum_{{i}=1}^{n}
(
x_{i} -( \beta_1 {a}_{i} + \beta_0 ))^2
\Big)
\end{align}
where
$(\hat{\beta}_0, \hat{\beta}_1)$
is called
"sample regression coefficients".

\end{itemize}
This is easily solved as follows.
Taking partial derivatives with respect to
$\beta_0$,
$\beta_1 $,
and
equating the results to zero,
gives the equations
(i.e.,
"normal equations"),
\begin{align}
&
\frac{\partial \sigma_{}^2
(\beta_1, \beta_2)}{\partial \beta_0}
=
{\sum_{{i}=1}^{n} ({}x_{i} - \beta_0 - \beta_1 {a}_{i}{})}=0,
\quad
({i}=1,...,{n})
\label{NE1}
\\
&
\frac{\partial \sigma_{}^2
(\beta_1, \beta_2)}{\partial \beta_1}
=
{\sum_{{i}=1}^{n} ({}x_{i} - \beta_0 - \beta_1 {a}_{i}{}){a}_{i}}=0,
\quad
({i}=1,...,{n})
\label{NE2}
\end{align}
Solving it, we get that
\begin{align}
&
\hat{\beta}_1
=
\frac{s_{{a}x}}{s_{{a}{a}}},
\quad
\hat{\beta}_0 =
\overline{x}-\frac{s_{{a}x}}{s_{{a}{a}}}\overline{{a}},
\quad
\hat{\sigma}^2
(
=
\frac{1}{{n}}
\sum_{{i}=1}^{n}
(
x_{i} -( \hat{\beta}_1 {a}_{i} + \hat{\beta}_0 ))^2
\Big)
=
s_{xx} - \frac{s_{{a}x}^2}{s_{{a}{a}}}
\label{S1}
\end{align}
where
\begin{align}
&
{\bar {a}}=\frac{{a}_1 + \cdots + {a}_{n}}{{n}}, \qquad
{\bar x}=\frac{x_1 + \cdots + x_{n}}{{n}},\qquad
\label{Notation1}
\\
&
s_{{a}{a}}=\frac{({a}_1 -{\bar {a}})^2 + \cdots + ({a}_{n} - {\bar {a}})^2}{{n}},
\quad
s_{xx}=\frac{(x_1 -{\bar x})^2 + \cdots + (x_{n} - {\bar x})^2}{{n}},
\quad
\label{Notation2}
\\
&
s_{{a}x}=\frac{({a}_1 -{\bar {a}})(x_1 -{\bar x}) + \cdots +
({a}_{n} - {\bar {a}})(x_{n} - {\bar x})}{{n}}.
\label{Notation3}
\end{align}
\par
\noindent
\par
\noindent
\bf
Remark 1
\rm
[Applied mathematics].
The above result is in applied mathematics and neither in statistics nor
in quantum language.
The purpose of this paper is to add
a quantum linguistic story to Problem 1
(i.e.,
the least-squared method)
in the framework of quantum language.
\par
\noindent

\par
\noindent
\subsection{
Quantum language
(Axioms
and
Interpretation)
}
\label{subsec:12QL}
\par
\noindent
\par
As mentioned in Remark 1, our purpose is to 
add a quantum linguistic story.
%
Thus,
we shall, according to ref.\cite{Ishi9}, mention the overview of quantum language (or, measurement theory, in short, MT).
\par
\par
\rm
Quantum language is characterized as the linguistic turn of the Copenhagen interpretation of quantum mechanics({\it cf.} refs. \cite{Ishi5},
{{{}}}{\cite{Neum}}).
Quantum language (or, measurement theory ) has two simple rules
(i.e. Axiom 1 (concerning measurement) and Axiom 2 ( 
concerning causal relation))
and the linguistic interpretation (= how to use the Axioms 1 and 2). 
That is,
\begin{itemize}
\item[(B$_1$)]$\underset{\mbox{(=MT(measurement theory))}}{\fbox{Quantum language}}
=
\underset{\mbox{(measurement)}}{\fbox{Axiom 1}}
+
\underset{\mbox{(causality)}}{\fbox{Axiom 2}}
+
\underset{\mbox{(how to use Axioms)}}{\fbox{linguistic interpretation}}
\label{eq1}
$
\end{itemize}
({\it cf.} refs.
{{{}}}{\cite{Ishi2}-\cite{Ishi10}}). This is all of quantum language.
\par
This theory is formulated in a certain $C^*$-algebra ${\cal A}$({\it cf.} ref.
{{{}}}{\cite{Saka}}), and is classified as follows:
\begin{itemize}
\item[(B$_2$)]
$
\quad
\underset{\text{\scriptsize }}{\text{Quantum language(=MT)}}
$
$\left\{\begin{array}{ll}
\text{quantum MT$\quad$(when ${\cal A}$ is non-commutative)}
\\
\\
\text{classical MT
$\quad$
(when ${\cal A}$ is commutative, i.e., ${\cal A}=C_0(\Omega)$)}
\end{array}\right.
$
\end{itemize}
where $C_0(\Omega)$
is
the $C^*$-algebra composed of all continuous 
complex-valued functions vanishing at infinity
on a locally compact Hausdorff space $\Omega$.

Since our concern in this paper is concentrated to 
classical systems,
we devote ourselves to the commutative $C^*$-algebra $C_0(\Omega)$,
which is quite elementary.
Therefore, we believe that all statisticians
can understand our assertion
(i.e.,
the quantum linguistic approach to statistics
).

Let $\Omega$ is a locally compact Hausdorff space, which is also called
a state space. And thus, an element $\omega (\in \Omega )$ is said to be a state.
Let $C(\Omega)$ be the $C^*$-algebra composed of all bounded continuous 
complex-valued functions on a locally compact Hausdorff space $\Omega$.
The norm $\| \cdot \|_{C(\Omega )}$ is usual, i.e.,
$\| f \|_{C(\Omega )} = \sup_{\omega \in \Omega } |f(\omega )|$
$(\forall f \in C(\Omega ))$.
\par
\rm
Motivated by Davies' idea ({\it cf.} ref.
{{{}}}{\cite{Davi}}) in quantum mechanics,
an observable ${\mathsf O}=(X, {\mathcal F}, F)$ in $C_0(\Omega )$
(or, precisely, in $C(\Omega )$) is defined as follows:
\begin{itemize}
\item[(C$_1$)]
$X$ is a topological space. 
${\mathcal F} ( \subseteq 2^X$(i.e., the power set of $X$) is a field,
that is, it satisfies the following conditions (i)--(iii):
(i):
$\emptyset \in {\cal F}$, 
(ii):$\Xi \in {\mathcal F} \Longrightarrow X\setminus \Xi  \in 
{\mathcal F}$,
(iii):
$\Xi_1, \Xi_2,\ldots, \Xi_{n} \in {\mathcal F} \Longrightarrow \cup_{k=1}^n \Xi_k \in {\mathcal F}$.
\item[(C$_2$)]
The map $F: {\cal F} \to C(\Omega )$ satisfies that 
\begin{align}
0 \le [F(\Xi )](\omega ) \le 1, \quad [F(X )](\omega )=1
\qquad
(\forall \omega \in \Omega )
\nonumber \end{align}
and moreover, 
if
\begin{align}
\Xi_1, \Xi_2,\ldots, \Xi_{k}, \ldots \in {\mathcal F},
\quad
\Xi_m \cap \Xi_n = \emptyset \quad( m \not= n ),
\quad
\Xi = \cup_{k=1}^\infty \Xi_k \in {\mathcal F},
\nonumber \end{align} 
then, it holds
\begin{align}
[F(\Xi)](\omega) = \lim_{n \to \infty } \sum_{k=1}^n [F(\Xi_k )](\omega )
\quad
(\forall \omega \in \Omega )
\nonumber \end{align}
\end{itemize}
Note that Hopf extension theorem
({\it cf.}
ref.
{{{}}}{\cite{Yosi}})
guarantees that
$(X, {\cal F}, [F(\cdot)](\omega))$
is regarded as the mathematical probability space.
\par

\vskip0.5cm

%
%
\par
\noindent
\bf
Example 1
\rm
[The normal observable].
Put $\Omega={\mathbb R} \times {\mathbb R}_+ =
\{ ( \mu, \sigma) \in {\mathbb R}^2 \; : \; \sigma > 0 \}$.
Define the normal observable
${\mathsf O}_G = ({\mathbb R}, {\mathcal B}_{\mathbb R}, {{{G}}} )$ 
in $C_0({\mathbb R} \times {\mathbb R}_+)$
such that
\par
\noindent
\begin{align}
&
[{{{G}}}
( \Xi)]
({}\omega{})
=
\frac{1}{{{\sqrt{2 \pi }\sigma{}}}}
\underset{{ \Xi }}{\int}
\exp[{}- \frac{ ({}{}{x} - {}{\mu}  {})^2 
}
{2 \sigma^2}    {}] d {}{x}
\label{eq2}
\\
&
\qquad 
({}\forall  \Xi \in {\cal B}_{{\mathbb R}}
\mbox{(=Borel field in ${\mathbb R}$)},
\quad
\forall   {}{\omega}=(\mu, \sigma )    \in \Omega = {\mathbb R}\times {\mathbb R}_+{}).
\nonumber
\end{align}
This observable is the most fundamental in this paper.
\vskip0.5cm

\par
Now we shall briefly explain "quantum language (B)" in classical systems as follows:
\par
A measurement of an observable
${\mathsf O}=(X, {\mathcal F}, F)$
for a system with a state $\omega (\in \Omega )$
is denoted by
${\mathsf M}_{C_0(\Omega)} ({\mathsf O}, S_{[\omega]})$.
By the measurement, a measured value $x (\in X)$ is obtained as follows:
\par
\noindent
\bf
Axiom 1
\rm
[Measurement].
\begin{itemize}
\item[(D$_1$)]
\sl
The probability that a measured value $x$
$( \in X)$ obtained by the measurement 
${\mathsf{M}}_{{{C_0(\Omega)}}} ({\mathsf{O}}$
${ \equiv} (X, {\cal F}, F),$
{}{$ S_{[\omega_0]})$}
belongs to a set 
$\Xi (\in {\cal F})$ is given by
$
[F(\Xi) ](\omega_0 )
$.
\end{itemize}
\rm
\par
\noindent
\par
\noindent
\bf
Axiom 2
\rm
[Causality].
\begin{itemize}
\item[(D$_2$)]
\sl
The causality is represented by a Markov operator
$\Phi_{21} : C_0(\Omega_2 ) \to C_0(\Omega_1 )$.
Particularly, the deterministic causality
is represented by a continuous map
$\phi_{12} : \Omega_1 \to \Omega_2$
such that
$$
f_2 ( \phi_{12}(\omega_1))
=
[\Phi_{12}(f_2)](\omega_1)
\quad
(\forall f_2 \in C_0(\Omega_2), \omega_1 \in \Omega_1 )
$$
\rm
Also, see (\ref{HeisPic}) later.
\end{itemize}
\par
\vskip0.3cm
\par
\noindent
\bf
Interpretation
\rm
[Linguistic interpretation].
Although there are several linguistic rules in quantum language, the following is the most important:
\begin{itemize}
\item[(D$_3$)]
\sl
Only one measurement is permitted.
And thus, the state is only one and does not move.
\end{itemize}
\rm
In order to read this paper,
it suffices to understand the above three.
For the further arguments, see refs.
{{{}}}{\cite{Ishi4},\cite{Ishi5},\cite{Ishi6},\cite{Ishi7},\cite{Ishi8},
\cite{Ishi9},
\cite{Ishi10}}.
\vskip0.5cm

\bf
\par
\noindent
Remark 2
\rm
[Random variables in Kolmogorov's probability theory].
It should be noted that
the word of
"random variable"
in Kolmogorov's probability theory
is not included in quantum language
(i.e.,
Axioms 1 and 2).
However,
the theory of random variables
(i.e.,
Kolmogorov's probability theory)
is frequently used in the mathematical proofs of quantum linguistic statements,
just like the mathematical theory of differential equations is used in the proofs of Newtonian mechanical statements.
(Continued to Remark 3).

\subsection{
Fisher's maximum likelihood method
(concerning Axiom 1)}
\label{subsec:13Fish}
\rm
\par
\noindent 
\par
It is usual to consider that
we do not know the pure state
$\omega_0$
$(
\in
\Omega
)$
when
we take a measurement
${\mathsf{M}}_{{{C_0(\Omega)}}} ({\mathsf{O}}, S_{[\omega_0]})$.
That is because
we usually take a measurement ${\mathsf{M}}_{{{C_0(\Omega)}}} ({\mathsf{O}},
S_{[\omega_0]})$
in order to know the state $\omega_0$.
Thus,
when we want to emphasize that
we do not know the state $\omega_0$,
${\mathsf{M}}_{{{C_0(\Omega)}}} ({\mathsf{O}}, S_{[\omega_0]})$
is denoted by
${\mathsf{M}}_{{{C_0(\Omega)}}} ({\mathsf{O}}, S_{[\ast]})$.

\vskip0.5cm
\par
\noindent
{\bf Theorem 1}
\rm
[Fisher's maximum likelihood method ({\it cf.} refs.
{{{}}}{\cite{Ishi3},\cite{Ishi4})].
\sl
Consider a measurement
${\mathsf M}_{{C_0(\Omega)}}
(
{\mathsf O}=(X , {\cal F} , F )
,$
$ S_{[*]}{})$.
Assume that
we know that the measured value $x \;(\in X )$
obtained by a measurement
${\mathsf M}_{{C_0(\Omega)}}
(
{\mathsf O}=(X , {\cal F} , F )
,$
$ S_{[*]}{})$
belongs to
$\Xi (\in {\cal F})$.
Then,
there is a reason to infer that the unknown state
$[\ast ]$ is equal to $\omega_0 (\in \Omega )$ such that
\begin{align}
\min_{\omega_1 \in \Omega}
\frac{[F(\Xi)](\omega_0)}
{
[F(\Xi)](\omega_1)
}
\Big(
=
\frac{[F(\Xi)](\omega_0)}
{
\max_{\omega_1 \in \Omega} [F(\Xi)](\omega_1)
}
\Big)
=
1
\label{eq12}
\end{align}
if the righthand side of this formula exists.
Also, if $\Xi=\{x\}$, it suffices to calculate the $\omega_0 (\in \Omega)$ 
such that
\begin{align}
L(x, \omega_0)
=1
\end{align}
where the likelihood function $L(x, \omega ) (\equiv L_x (\omega ))$ is defined by
\begin{align}
L(x, \omega )
=
\inf_{\omega_1 \in \Omega }
\Big[\lim_{\Xi \supseteq \{x \}, \; [F(\Xi)](\omega_1) \not= 0,\; \Xi \to \{x \}} \frac{[F(\Xi)](\omega)}{
[F(\Xi)](\omega_1)}
\Big]
\label{eq13}
\end{align}
\par
\noindent
\rm
\par
\noindent

\par
\noindent
\bf
Definition 1
\rm
[Product observable (or, simultaneous observable),
simultaneous measutement].
For each $k=1,2,\cdots, K$,
consider an observable
${\mathsf{O}_k}$
$=(X_k ,$
${\cal F}_k,$
$ 
{F}_k)$
in ${C_0(\Omega)}$.
Define the simultaneous observable
$\bigtimes_{k=1}^K {\mathsf{O}_k}$
$=(\bigtimes_{k=1}^K X_k ,$
$ \boxtimes_{k=1}^K {\cal F}_k,$
$ 
\bigtimes_{k=1}^K {F}_k)$
in ${C_0(\Omega)}$ such that
\begin{eqnarray}
(\bigtimes_{k=1}^K {F}_k)(\bigtimes_{k=1}^K {\Xi}_k )
=
F_1(\Xi_1) F_2(\Xi_2) \cdots F_K(\Xi_K)
\label{simu}
\\
\;
(
\forall \Xi_k \in {\cal F}_k,
\forall k=1,...,K
).
\qquad
\qquad
\nonumber
\end{eqnarray}
where
$ \boxtimes_{k=1}^K {\cal F}_k$ is the product field of ${\cal F}_k$
$(k=1,2, \cdots, K)$.

\par
\vskip1.0cm
For each
$k=1,$
$2,...,K$,
consider a measurement
${\mathsf{M}}_{{{C_0(\Omega)}}} ({\mathsf{O}_k}$
${\; :=} (X_k, {\cal F}_k, F_k),$
$ S_{[\omega]})$.
However,
since
the linguistic interpretation (D$_3$)
says that
only one measurement is permitted,
the
multiple measurements
$\{
{\mathsf{M}}_{{{C_0(\Omega)}}} ({\mathsf{O}_k},S_{[\omega]})
\}_{k=1}^K$
are prohibited.
Thus,
this
$\{
{\mathsf{M}}_{{{C_0(\Omega)}}} ({\mathsf{O}_k},S_{[\omega]})
\}_{k=1}^K$
is
represented by the simultaneous measurement
${\mathsf{M}}_{{{{C_0(\Omega)}}}} (
\bigtimes_{k=1}^K {\mathsf{O}_k}$,
$ S_{[\omega]})$.

\par
\noindent
\bf
Example 2
\rm
[Simultaneous normal observable].
Let
${\mathsf O}_G = ({\mathbb R}, {\mathcal B}_{\mathbb R}, {{{G}}} )$ be the normal observable
in $C_0({\mathbb R} \times {\mathbb R}_+)$ in Example 1.
Let $n$ be a natural number.
Then, we obtain the simultaneous normal observable
${\mathsf O}_G^n = ({\mathbb R}^n, {\mathcal B}_{\mathbb R}^n, {{{G}}^n} )$ 
in $C_0({\mathbb R} \times {\mathbb R}_+)$.
That is,
\par
\noindent
\begin{align}
&
[{{{G}}}^n
(\bigtimes_{{i}=1}^n \Xi_{i})]
({}\omega{})
=
[{{{G}}}^n
(\bigtimes_{{i}=1}^n \Xi_{i})]
({}\mu, \sigma{})
=
\frac{1}{({{\sqrt{2 \pi }\sigma{}}})^n}
\underset{{\bigtimes_{{i}=1}^n \Xi_{i} }}{\int \cdots \int}
\exp[{}- \frac{\sum_{{i}=1}^n ({}{}{x_{i}} - {}{\mu}  {})^2 
}
{2 \sigma^2}    {}] d {}{x_1} d {}{x_2}\cdots dx_n
\label{eq2}
\\
&
\qquad 
({}\forall  \Xi_{i} \in {\cal B}_{{\mathbb R}}
\mbox{(=Borel field in ${\mathbb R}$)},
({}{i}=1,2,\ldots, n),
\quad
\forall   {}{\omega}=(\mu, \sigma )    \in \Omega = {\mathbb R}\times {\mathbb R}_+{}).
\nonumber
\end{align}
Thus, we have the simultaneous measurement
${\mathsf M}_{C_0({\mathbb R} \times {\mathbb R}_+ )} ({\mathsf O}_G^n = ({\mathbb R}^n, {\mathcal B}_{\mathbb R}^n, {{{G}}^n} )$,
$S_{[\ast]}{})$ 
in $C_0({\mathbb R} \times {\mathbb R}_+)$.
Assume that a measured value $x=(x_1, x_2, \ldots, x_n ) (\in
{\mathbb R}^n )$ is obtained by the measurement.
Since the likelihood function
$L_x(\mu, \sigma)(=L(x, (\mu,\sigma)) $
is defined by
\begin{align}
L_x(\mu, \sigma)
&
=
\frac{C}{({{\sqrt{2 \pi }\sigma{}}})^n}
\exp[{}- \frac{\sum_{{i}=1}^n ({}{}{x_{i}} - {}{\mu}  {})^2 
}
{2 \sigma^2}    {}] 
\quad
(\mbox{the constant $C$ is independent of $\mu, \sigma$})
\end{align}
it suffices to calculate the following equations:
\begin{align}
\frac{\partial L_x(\mu, \sigma)}{\partial \mu}=0,
\quad
\frac{\partial L_x(\mu, \sigma)}{\partial \sigma}=0
\label{eq15}
\end{align}
Thus, we see, by
Theorem 1
(
Fisher's maximum likelihood method),
that
the unknown state
$[\ast]
$
can be inferred by
$(\hat{\mu}, \hat{\sigma})$,
that is,
%
\begin{align}
&
\hat{\mu}
(x) =
\hat{\mu}
(x_1,x_2,\ldots , x_n ) =
\frac{x_1 + x_2 + \cdots + x_n}{n}
\label{eq3}
\\
&
{{\hat{\sigma}}}
(x) =
{{\hat{\sigma}}}
(x_1,x_2,\ldots , x_n ) =
\sqrt{
\frac
{\sum_{{i}=1}^n ( x_{i} - 
\hat{\mu}
(x))^2}
{n}
}
\label{eq4}
\end{align}
For example, consider the following
image observable:
$\hat{\mu}({\mathsf O}_G^n) $
$= ({\mathbb R}, {\mathcal B}_{\mathbb R}, {{{G}}^n} \circ \hat{\mu}^{-1} )$
in $C_0({\mathbb R} \times {\mathbb R}_+)$
such that
\begin{align}
[({{{G}}^n} \circ \hat{\mu}^{-1})(\Xi_1)](\omega)
=
&
\frac{1}{({{\sqrt{2 \pi }\sigma{}}})^n}
\underset{
\{ x \in {\mathbb R}^n \;:\; {\hat{\mu}}(x) \in \Xi_1 \}}
{\int \cdots \int}
\exp[{}- \frac{\sum_{{i}=1}^n ({}{}{x_{i}} - {}{\mu}  {})^2 
}
{2 \sigma^2}    {}] d {}{x_1} d {}{x_2}\cdots dx_n
\label{Gauss1}
\intertext{which is calculated as follows:}
=
&
\frac{\sqrt{n}}{{\sqrt{2 \pi }\sigma{}}}
\int_{{\Xi_1}} \exp[{}- \frac{n({}{}{x} - {}{\mu}  {})^2 }{2 \sigma^2}    {}] d {}{x}
\label{Gauss2}
\\
&
\quad
(
\forall   {}{\omega} =(\mu, \sigma)   \in \Omega \equiv {\mathbb R}{}\times {\mathbb R}_+,
\forall
\Xi_1 \in {\mathcal B}_{\mathbb R}).
\nonumber
\end{align}
%

\bf
\par
\noindent
Remark 3
\rm
[Kolmogorov's probability theory].
Although the derivation of (\ref{Gauss2}) from (\ref{Gauss1})
may not be easy,
it is the problem in mathematics.
Although there are several derivations,
the calculation in the framework of Kolmogorov's probability theory (ref. \cite{Kolm}) may be the most elegant.
Thus, mathematical theories
(e.g.,
Kolmogorov's probability theory,
operator theory (ref. \cite{Saka})
) are frequently used in quantum language.
%

\subsection{The Heisenberg picture (concerning Axiom 2)}
\label{subsec;14Heisenberg}

\par
\noindent
\par
Consider a tree-like ordered set
$(T{\; :=}\{t_0, t_1, ..., t_n \},$
$ \le )$
with the root $t_0$
(i.e.,
$t_0 \le t \;(\;\forall t \in T)$).
This is also characterized by
the parent map
$\tau: T \setminus \{t_0\} \to T$
such that
$\tau( t)= \max \{ s \in T \;|\; s < t \}$.
Put
$T_\le^2=\{
(t,t') \in T^2 \;:\;
t \le t' \}$.
In Figure \ref{Figure1},
see
the root
$t_0$,
the parent map:
$\tau(t_3)=\tau(t_4)=t_2$,
$\tau(t_2)=\tau(t_5)=t_1$,
$\tau(t_1)=\tau(t_6)=\tau(t_7)=t_0$
\par
\noindent
\begin{figure}[htbp]
\setlength{\unitlength}{0.7mm}
\begin{picture}(100,60)(40,0)
\put(60,0){
\put(33,18){\makebox(10,10)[r]{${t_0}$}}
\put(65,30){\makebox(10,10)[r]{${t_1}$}}
\put(93,45){\makebox(10,10)[r]{${t_2}$}}
\put(131,50){\makebox(10,10)[r]{${t_3}$}}
\put(131,35){\makebox(10,10)[r]{${t_4}$}}
\put(91,14){\makebox(10,10)[r]{${t_5}$}}
\put(60,10){\makebox(10,10)[r]{${t_6}$}}
\put(60,0){\makebox(10,10)[r]{${t_7}$}}
\put(60,33){\vector(-3,-1){15}} 
\put(60,20){\vector(-3,1){15}} 
\put(60,7){\vector(-3,2){15}} 
\put(90,48){\vector(-3,-2){13}} 
\put(90,20){\vector(-3,2){13}} 
\put(128,55){\vector(-3,-1){13}} 
\put(128,40){\vector(-3,2){13}} 
\put(50,16){$\tau$}
\put(47,32){$\tau$}
\put(47,5){$\tau$}
\put(73,45){$\tau$}
\put(83,25){$\tau$}
\put(118,57){$\tau$}
\put(118,38){$\tau$}
}
\end{picture}
\caption{
\label{Figure1}
Tree
}
\end{figure}

For each $t \in T$, a commutative $C^*$-algebra $C_0(\Omega_t )$ is associated.
According to Axiom 2,
consider a Markov relation (i.e., causal relation)
$\{ \Phi_{t, t'} : C_0(\Omega_{t'})  \to C_0({\Omega}_{t})  \}_{ (t,t')\in
T_\le^2}$, which is also represented by
$\{ \Phi_{\tau(t), t} : C_0(\Omega_{t})  \to C_0(\Omega_{\tau(t)})  \}_{ 
t \in T \setminus \{t_0\}}$.
In this paper, we consider the deterministic case,
that is, the case that $\Phi_{\tau(t), t} : C_0(\Omega_{t})  \to C_0(\Omega_{\tau(t)})$ is represented by the continuous map
(called "{\it causal map}")
$\phi_{\tau(t), t} : \Omega_{\tau(t)} \to \Omega_{t}$ such as
\begin{align}
(\Phi_{\tau(t), t} f_t)(\omega_{\tau(t)})
=
f_t ( \phi_{\tau(t), t}(\omega_{\tau(t)}))
\quad
(\forall \omega_{\tau(t)} \in \Omega_{{\tau(t)}},
\forall
f_t \in C_0(\Omega_t ), \forall t \in T \setminus \{t_0\})
\label{HeisPic}
\end{align} 
Let an observable ${\mathsf O}_t{\; :=}
(X_t, {\cal F}_{t}, F_t)$ in the $C_0(\Omega_t)$ 
be given for each $t \in T$.
$\Phi_{\tau(t), t} {\mathsf O}_t$
is defined by
$(X_t,$
$ {\cal F}_{t},$
$ \Phi_{\tau(t), t}F_t)$ in the $C_0(\Omega_{\tau(t)})$.
And let $\omega_{0} \in \Omega_{t_0}$.
Consider {\lq\lq}measurements" such as
\begin{itemize}
\item[(E)]
for each $t \in T$, take a measurement of an observable
${\mathsf O}_t$ for the system with a {\lq\lq}moving state"
$\phi_{t_0, t}(\omega_{0}) \in \Omega_{t}$.
\end{itemize}
where the meaning of {\lq\lq}moving state" is not clear yet.
Recalling that
the linguistic interpretation (D$_3$) says that a state never moves,
we consider the meaning of the (E) as follows:
For each $s \in T$,
put $T_s =\{ t \in T \;|\; t \ge s\}$.
And define the observable
${\widehat{\mathsf O}}_s
=(\bigtimes_{t \in T_s}X_t, \boxtimes_{t \in T_s}{\cal F}_t, {\widehat{F}}_s)$
in $C_0(\Omega_s)$
(due to the Heisenberg picture)
as follows:
\par
\noindent
\begin{align}
\widehat{\mathsf O}_s
&=
\left\{\begin{array}{ll}
{\mathsf O}_s
\quad
&
\!\!\!\!\!\!\!\!\!\!\!\!\!\!\!\!\!\!
\text{(if $s \in T \setminus \tau (T) \;${})}
\\
\\
{\mathsf O}_s
\bigtimes
({}\bigtimes_{t \in \tau^{-1} ({}\{ s \}{})} \Phi_{ \tau(t), t} \widehat {\mathsf O}_t{})
\quad
&
\!\!\!\!\!\!
\text{(if $ s \in \tau (T) ${})}
\end{array}\right.
\label{prod}
\end{align}
Using (\ref{prod}) iteratively,
we can finally obtain the observable
$\widehat{\mathsf O}_{t_0}$
in $C_0(\Omega_{t_0})$.
Thus, the above (E) is represented by
the measurement
${\mathsf{M}}_{{C_0(\Omega_{t_0} )}} (\widehat{\mathsf{O}}_{t_0},$
$ S_{[\omega_0]})$.
Since the causal map is assumed to be deterministic in this paper, 
the $\widehat{\mathsf O}_{t_0}$
is simply represented by the simultaneous observable such as
$\widehat{\mathsf O}_{t_0}$
$=\bigtimes_{t \in T} \Psi_{t_0, t} {\mathsf O}_{t}$
({\it cf.} refs.
{{{}}}{\cite{Ishi4},\cite{Ishi7},\cite{Ishi8}).
\par
\vskip0.5cm

\par
\noindent
\bf
Remark 4
\rm
[What is regression analysis?].
Since regression analysis has various aspects,
it is not easy to answer the question:
"What is regression analysis?"
However,
we can say that
regression analysis is at least related to
the inference concerning
${\mathsf{M}}_{{C_0(\Omega_{t_0} )}} (\widehat{\mathsf{O}}_{t_0},$
$ S_{[\ast]})$.
In this sense,
regression analysis must be related to
Axiom 2 as well as Axiom 1.
On the other hand, Fisher's maximum likelihood method is related to only Axiom 1.
We believe that
the reason that
regression analysis is famous is to be related to Axiom 2. 
As seen in (B$_1$),
the importance of Axiom 2 (Causality) is explicitly emphasized in quantum language
and not in statistics.
Thus, we think that
regression analysis plays the role of Axiom 2
in the conventional statistics.
And,
in Section 2,
we will point out that
the  term "explanatory variable"
is understood as a kind of causal map
in quantum language.

\par
\vskip1.0cm
\par
\rm
\subsection{The reverse relation between confidence interval and statistical hypothesis testing}
\label{subsec:15Hypo}
\par
\noindent
\par
Let
${\mathsf O} = ({}X, {\cal F} , F{}){}$
be an observable
formulated in a
commutative $C^*$-algebra
${C_0(\Omega)}$.
Let $X$ be  a topological space.
Let $\Theta$ be a locally compact space with the 
semi-distance $d^x_{\Theta}$
$(\forall x \in X)$,
that is,
for each $x\in X$,
the map
$d^x_{\Theta}: \Theta^2 \to [0,\infty)$
satisfies that
(i):$d^x_\Theta (\theta, \theta )=0$,
(ii):$d^x_\Theta (\theta_1, \theta_2 )$
$=d^x_\Theta (\theta_2, \theta_1 )$,
(ii):$d^x_\Theta (\theta_1, \theta_3 )$
$\le d^x_\Theta (\theta_1, \theta_2 )
+
d^x_\Theta (\theta_2, \theta_3 )
$.

Let
${\widehat{E}}:X \to \Theta$
and
$\pi: \Omega \to \Theta$ be
continuous maps,
which are respectively called
an estimator and
a quantity.
Let
$\alpha$
be a real number such that
$0 < \alpha \ll 1$,
for example,
$\alpha = 0.05$.
For any state
$ \omega ({}\in \Omega)$,
define
the positive number
$\eta^\alpha_{\omega}$
$({}> 0)$
such that:
\begin{align}
\eta^\alpha_{\omega}
&
=
\inf
\{
\eta > 0:
[F(\{ x \in X \;:\; 
d^x_\Theta ( {\widehat{E}}(x) , \pi( \omega ) )
\ge \eta
\}
)](\omega )
\le \alpha
\}
\label{Defeta}
\\
\Big(
&=
\inf
\{
\eta > 0:
[F(\{ x \in X \;:\; 
d^x_\Theta ( {\widehat{E}}(x) , \pi( \omega ) )
< \eta
\}
)](\omega )
\ge 1- \alpha
\}
\Big)
\nonumber
\end{align}
Then Axiom 1 says that:
\rm
\begin{enumerate}
\item[(F$_1$)]
\it
\sl
the probability,
that
the measured value $x$
obtained
by the measurement
${\mathsf M}_{C_0(\Omega)} \big({}{\mathsf O}:= ({}X, {\cal F} , F{})  ,$
$ S_{[\omega_0 {}] } \big)$
satisfies the following
condition (\ref{eq10}),
is more than or equal to
$1-\alpha$
({}e.g., $1-\alpha= 0.95${}).
\begin{align}
d^x_\Theta ({\widehat{E}}(x),  \pi(\omega_0){}) <  {\eta }^\alpha_{\omega_0}
\label{eq10}
\end{align}
\end{enumerate}
or equivalently,
\begin{enumerate}
\item[(F$_2$)]
\it
\sl
the probability,
that
the measured value $x$
obtained
by the measurement
${\mathsf M}_{C_0(\Omega)} \big({}{\mathsf O}:= ({}X, {\cal F} , F{})  ,$
$ S_{[\omega_0 {}] } \big)$
satisfies the following
condition (\ref{eq11}),
is less than or equal to
$\alpha$
({}e.g., $\alpha= 0.05${}).
\begin{align}
d^x_\Theta ({\widehat{E}}(x),  \pi(\omega_0){}) \ge  {\eta }^\alpha_{\omega_0}
\label{eq11}
\end{align}
\end{enumerate}

\rm
\par
\noindent
\bf
Theorem 2
\rm
[Confidence interval and statistical hypothesis testing
({\it cf.} ref. \cite{Ishi9})
].
\sl
Let
${\mathsf O} = ({}X, {\cal F} , F{}){}$
be an observable
formulated in a
commutative $C^*$-algebra
${C_0(\Omega)}$.
Let
${\widehat{E}}:X \to \Theta$
and
$\pi: \Omega \to \Theta$ be
an estimator and
a quantity
respectively.
Let $\eta_\omega^\alpha$ be as defined in the formula (\ref{Defeta}).
\par
\noindent
From the ($F_1$), we assert "the confidence interval method" as follows:
\begin{enumerate}
\rm
\item[(G$_1$)]
[The confidence interval method].
\sl
\it
\sl
For any $x \in X$, define
\begin{align}
I_{x}^{1- \alpha}
=
\{
\pi(\omega)
(\in
\Theta)
:
d^x_\Theta ({}{\widehat{E}}(x),
\pi(\omega )
)
<
\eta^{1- \alpha}_{\omega }
\}
\label{eq12} 
\end{align}
which is called
\it
the $({}1- \alpha{})$-confidence interval. 
\sl
Let
$x (\in X)$
be
a measured value $x$
obtained
by the measurement
${\mathsf M}_{C_0(\Omega)} \big({}{\mathsf O}:= ({}X, {\cal F} , F{})  ,$
$ S_{[\omega_0 {}] } \big)$.
Then, the probability
that
$I_x^{1-\alpha} \ni \pi(\omega_0)$
is more than or equal to $1- \alpha$.
\end{enumerate}
From the ($F_2$), we assert "the statistical hypothesis test" as follows:
\rm
\begin{itemize}
\item[(G$_2$)]
[The statistical hypothesis test].
\sl
Assume that
a state $\omega_0$
satisfies that
$\pi(\omega_0)
\in
H_N
( \subseteq \Theta )
$,
where $H_N$
is called a 
"null hypothesis".
Put
\begin{align}
&
{\widehat R}_{{H_N}}^{\alpha; \Theta}
=
\bigcap_{\omega \in  \Omega \mbox{ \footnotesize such that }
\pi(\omega) \in {H_N}}
\{
{\widehat{E}}({x})
(\in
\Theta)
:
d^x_\Theta ({}{\widehat{E}}(x),
\pi(\omega )
)
\ge
\eta^\alpha_{\omega }
\}
\label{eq13}
\intertext{and also}
&
{\widehat R}_{{H_N}}^{\alpha; X}
=
{\widehat{E}}^{-1}(
{\widehat R}_{{H_N}}^{\alpha; \Theta})
=
\bigcap_{\omega \in  \Omega \mbox{ \footnotesize such that }
\pi(\omega) \in {H_N}}
\{
x
(\in
X)
:
d^x_\Theta ({}{\widehat{E}}(x),
\pi(\omega )
)
\ge
\eta^\alpha_{\omega }
\}
\label{eq114}
\end{align}
which is respectively called
\it
the $({}\alpha{})$-rejection region
of
the null hypothesis
${H_N}$.
\sl
Then,
the probability,
that
the measured value $x (\in X)$
obtained
by the measurement
${\mathsf M}_{C_0(\Omega)} \big({}{\mathsf O}:= ({}X, {\cal F} , F{}),$
%
$S_{[\omega_0]} \big)$
$($
where it should be noted that $\pi(\omega_0) \in H_N )$
satisfies the following
condition (\ref{eq15}),
is less than or equal to
$\alpha$
({}e.g., $\alpha= 0.05${}).
\begin{align}
"{\widehat{E}}(x) \in
{\widehat R}_{{H_N}}^{\alpha; \Theta}"
\mbox{  or equivalently  }
"
x
\in
{\widehat R}_{{H_N}}^{\alpha; X}"
\label{eq15}
\end{align}
\end{itemize}

\par
\vskip1.0cm
\par

\section{Regression analysis in quantum language}
\par
\noindent
\par

\def\Cal{\cal}
\def\bigstimes{\text{\large $\: \boxtimes \,$}}

In this section,
we show that
the least squared method
(mentioned in Section \ref{subsec:11Least}) acquires 
a quantum linguistic story as follows.
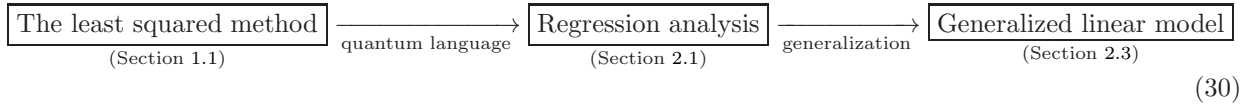
\begin{align}
\underset{\mbox{\scriptsize{(Section \ref{subsec:11Least})}}}{\fbox{The least squared method}}
\xrightarrow[\mbox{\scriptsize{quantum language}}]{}
\underset{\mbox{\scriptsize{(Section \ref{subsec:21GLM})}}}{\fbox{Regression analysis}}
\xrightarrow[\mbox{\scriptsize{generalization}}]{}
\underset{\mbox{\scriptsize{(Section \ref{subsec:23GLM})}}}{\fbox{Generalized linear model}}
\end{align}
Note that Theorem 1 (Fisher's maximum likelihood method)
and
Theorem 2
(Confidence interval and hypothesis test)
are only related to
Axiom 1.
On the other hand,
it should be noted that
Axiom 2
(as well as Axiom 1)
is used in regression analysis.


\subsection{Simple regression analysis in quantum language}
\label{subsec:21GLM}

\par
\noindent



\par
\noindent
\par
Put
$T=\{ 0,1,2, \cdots, i , \cdots, n \}$.
And let $(T, \tau: T \setminus \{0\} \to T )$ be the tree-like ordered set
(with the parallel structure)
such that
\begin{align}
\tau(i)=0
\qquad
(\forall i =1,2, \cdots, n)
\label{Tree}
\end{align}

\par
\noindent
For each $i \in T$,
define a locally compact space $\Omega_i$ such that
\begin{align}
&
\Omega_{0}={\mathbb R}^{2}
=
\Big\{
\beta
=\begin{bmatrix}
\beta_0 \\
\beta_1 \\
\end{bmatrix}
\;:\;
\beta_0, \beta_1 \in {\mathbb R}
\Big\}
\quad
\\
&
\Omega_{i}={\mathbb R}
=
\Big\{
\mu_i
\;:\;
\mu_i \in {\mathbb R}
\Big\}
\quad(i=1,2, \cdots, n )
\end{align}
Assume that
\begin{align}
a_i \in {\mathbb R} 
\qquad
(i=1,2, \cdots, n ),
\label{Explanatory1}
\end{align}
which are called {\it explanatory variables}
in the conventional statistics.
Consider the deterministic causal map
$\psi_{a_i}: \Omega_0(={\mathbb R}^2)
\to
\Omega_{i} (={\mathbb R})$
such that
\begin{align}
\Omega_0={\mathbb R}^2 \ni \beta =(\beta_0, \beta_1 ) \mapsto
\psi_{a_i} ( \beta_0, \beta_1)
=
\beta_0 + \beta_1 a_i= \mu_i \in \Omega_i ={\mathbb R}
\end{align}
which is equivalent to the deterministic Markov operator
$\Psi_{a_i}:
C_0(\Omega_{i}) \to C_0(\Omega_0)$
such that
\begin{align}
[{\Psi_{a_i}}(f_i)](\omega_0)
=
f_i( \psi_{a_i} (\omega_0))
\quad
(\forall f_i \in C_0(\Omega_{i}),
\;\;
\forall \omega_0 \in \Omega_0,
\forall i \in 1,2, \cdots, n)
\end{align}
Thus, under the identification: $a_i \Leftrightarrow \Psi_{a_i}$,
the term "explanatory variable" means a kind of causal relation $\Psi_{a_i}$.

\par
\noindent
\begin{figure}[htbp]
{
\setlength{\unitlength}{0.7mm}
\begin{picture}(80,55)(0,20)
\put(40,0){
\put(83,60){\makebox(10,10)[r]{${C_0 (\Omega_1(\equiv{\mathbb R}))}$}}
\put(83,45){\makebox(10,10)[r]{${C_0 (\Omega_2(\equiv{\mathbb R}))}$}}
\put(83,22){\makebox(10,10)[r]{${C_0 (\Omega_n(\equiv{\mathbb R}))}$}}
\put(27,41){\makebox(10,10)[r]{${C_0 (\Omega_0(\equiv{\mathbb R}^2))}$}}
\put(60,65){\vector(-3,-2){15}}
\put(60,50){\vector(-3,-1){13}}
\put(60,27){\vector(-3,2){15}}
%
\put(47,67){$\Psi_{a_1}$}
\put(50,51){$\Psi_{a_2}$}
\put(50,43){$\cdots \cdots$}
\put(50,38){$\cdots \cdots$}
\put(47,25){$\Psi_{a_n}$}
}
\end{picture}
}
\caption{
\label{Fig2}
Parallel structure
(Causal relation $\Psi_{a_i}$)
}
\end{figure}

\par
\noindent
For each
$i=1,2, \cdots, n$, define
the {\it
normal observable}
${\mathsf O}_{i} {{\equiv}} ({\mathbb R}, {\cal B}_{{\mathbb R}}, G_{\sigma_{}})$
in
$C_0 (\Omega_{i}(\equiv
{\mathbb R}))$
such that
\begin{align}
&
[G_{\sigma_{}}(\Xi)] (\mu ) = \frac{1}{(\sqrt{2 \pi \sigma_{}^2})}
\underset{\Xi}{\int} \exp
\Big[{- \frac{
(x -\mu )^2
}{2 \sigma_{}^2}}
\Big] dx
\qquad
(\forall \Xi \in {\cal B}_{{\mathbb R}}, \forall 
\mu
\in \Omega_{i} (\equiv {\mathbb R} ))
\end{align}
where $\sigma $ is a positive constant.
\par
\noindent
Thus, we have
the observable
${\mathsf O}_{0}^{a_i} {{\equiv}} ({\mathbb R}, {\cal B}_{{\mathbb R}}, \Psi_{a_i}G_{\sigma_{}})$
in
$C_0 (\Omega_{0}(\equiv
{\mathbb R}^2))$
such that
\begin{align}
&
[\Psi_{a_i}(G_{\sigma_{}}(\Xi))] (\beta ) = 
[(G_{\sigma_{}}(\Xi))] (\psi_{a_i}(\beta ))
=
\frac{1}{(\sqrt{2 \pi \sigma_{}^2})}
\underset{\Xi}{\int} \exp
\Big[{- \frac{   (x - (\beta_0 + a_{i{}} \beta_1 ))^2}{2 \sigma_{}^2}}
\Big] dx
\label{54}
\\
&
\qquad
(\forall \Xi \in {\cal B}_{{\mathbb R}}, \forall 
\beta
=(\beta_0, \beta_1  )\in \Omega_{0} 
(\equiv {\mathbb R}^{2} )
\nonumber
\end{align}
Hence, we have the simultaneous observable
$\bigtimes_{i=1}^n{\mathsf O}_{0}^{a_i} {{\equiv}} ({\mathbb R}^n, {\cal B}_{{\mathbb R}^n}, \bigtimes_{i=1}^n \Psi_{a_i}G_{\sigma_{}})$
in
$C_0 (\Omega_{0}(\equiv
{\mathbb R}^2))$
%
such that
\begin{align}
&
[(\bigtimes_{i=1}^n \Psi_{a_i}G_{\sigma_{}})
(\bigtimes_{i=1}^n \Xi_i)](\beta)
=
\bigtimes_{i=1}^n \Big(
[\Psi_{a_i}G_{\sigma_{}})
(\Xi_i)](\beta)\Big)
\nonumber
\\
=
&
\frac{1}{(\sqrt{2 \pi \sigma_{}^2})^n}
\underset{\bigtimes_{i=1}^n  \Xi_i}{\int \cdots \int} \exp
\Big[{- \frac{ \sum_{i=1}^n  (x_i - (\beta_0 + a_{i{}} \beta_1 ))^2}{2 \sigma_{}^2}}
\Big] dx_1 \cdots dx_n
\nonumber
\\
=
&
\underset{\bigtimes_{i=1}^n  \Xi_i}{\int \cdots \int} 
p_{(\beta_0, \beta_1, \sigma )}
(x_1, x_2, \cdots, x_n )
dx_1 \cdots dx_n
\label{Fundamental22}
\\
&
\qquad \qquad
\qquad
(\forall \bigtimes_{i=1}^n  \Xi_i \in {\cal B}_{{\mathbb R}^n}, 
\forall 
\beta
=(\beta_0, \beta_1 ) \in \Omega_{0} (\equiv {\mathbb R}^{2} )
)
\nonumber 
\end{align}
Assuming that $\sigma$ is variable, we have the observable
${\mathsf O}=
\Big({\mathbb R}^n(=X) , {\mathcal B}_{{\mathbb R}^n}(={\mathcal F}),
F \Big)$ in $C_0 ( \Omega_0 \times {\mathbb R}_+ )$
such that
\begin{align}
[F(\bigtimes_{i=1}^n \Xi_i )](\beta, \sigma_{})
=
[(\bigtimes_{i=1}^n \Psi_{a_i}G_{\sigma_{}})
(\bigtimes_{i=1}^n \Xi_i)](\beta)
\quad
(\forall \Xi_i \in {\cal B}_{{\mathbb R}},
\forall
(\beta , \sigma_{} ) \in {\mathbb R}^2(\equiv \Omega_0) \times {\mathbb R}_+ )
\end{align}

\par

\par
\noindent
\bf
Problem 2
\rm
[Simple regression analysis in quantum language]
\sl
Assume that
a measured value 
$x=\begin{bmatrix}
x_1
\\
x_2
\\
\vdots
\\
x_n
\end{bmatrix}
\in X={\mathbb R}^n $
is obtained by the measurement
${\mathsf M}_{C_0(\Omega_{0} \times {\mathbb R}_+)}(
{\mathsf O} \equiv (X, {\cal F}, F)
, S_{[(\beta_0,\beta_1, \sigma)]}
{}
)$.
We do not know the state
$
(\beta_0, \beta_1, \sigma_{}^2
)$.
%
%
%
Then, from the measured value
$
x=(
x_1, x_2, \ldots, x_n )
\in {\mathbb R}^{n}$,
infer the $\beta_0, \beta_1, \sigma_{}$!
That is,
represent the
$(\beta_0, \beta_1, \sigma_{})$
by $(\hat{\beta}_0(x), \hat{\beta}_1(x), \hat{\sigma}_{}(x))$
(i.e.,
the functions of $x$).
\rm
\vskip0.5cm
\par
\noindent
\bf
Answer.
\rm
%
%
%
Taking partial derivatives with respect to
$\beta_0$,
$\beta_1 $,
$\sigma_{}^2$,
and
equating the results to zero,
gives the $\log$-likelihood equations.
That is, putting
$L(\beta_0, \beta_1, \sigma^2, x_1, x_2, \cdots, x_n)=\log p_{(\beta_0, \beta_1, \sigma )}
(x_1, x_2, \cdots, x_n )$,
we see that
\begin{align}
&
\frac{\partial L}{\partial \beta_0}=0 \quad \Longrightarrow
\quad
{\sum_{i=1}^n 
{(x_i - (\beta_0 + {{}} a_{i{}} \beta_1 ))}
}
=0
\label{NNNE1}
\\
&
\frac{\partial L}{\partial \beta_1}=0 \quad
\Longrightarrow
\quad
{\sum_{i=1}^n 
{a_{i{}}(x_i - (\beta_0 + {{}} a_{i{}} \beta_1 ))}
}
=0
\label{NNNE2}
\\
&
\frac{\partial L}{\partial \sigma^2}=0 \quad
\Longrightarrow
-\frac{n}{2\sigma_{}^2}
+
\frac{1}{2\sigma_{}^4}
{\sum_{i=1}^n ({}x_i - \beta_0 - \beta_1 a_i{})^2 }
=0
\label{NNNE3}
\end{align}
\par
\par
\noindent
Therefore,
using the notations
(\ref{Notation1})-(\ref{Notation3}),
we obtain that
%
\begin{align}
&
\hat{\beta}_0(x)=\overline{x}
-\hat{\beta}_1(x) \overline{a}=\overline{x}
-\frac{s_{ax}}{s_{aa}} \overline{a},
\quad
\hat{\beta}_1(x)=\frac{s_{ax}}{s_{aa}}
\label{SSSS1}
\intertext{and}
&
(\hat{\sigma_{}}(x))^2=
\frac{\sum_{i=1}^n \Big( x_i - (
\hat{\beta}_0 (x)+ a_{i{}} \hat{\beta}_1 (x) )
\Big)^2}{n}
\nonumber
\\
=
&
\frac{\sum_{i=1}^n \Big( x_i - (
\overline{x}
-\frac{s_{ax}}{s_{aa}} \overline{a})
-
a_i
\frac{s_{ax}}{s_{aa}} 
\Big)^2}{n}
=
\frac{\sum_{i=1}^n \Big(( x_i - 
\overline{x})
+(
\overline{a}
-
a_i)
\frac{s_{ax}}{s_{aa}} 
\Big)^2}{n}
\nonumber
\\
=
&
s_{xx}
-
2 s_{ax}\frac{s_{ax}}{s_{aa}} 
+
s_{aa}(\frac{s_{ax}}{s_{aa}})^2 
=
s_{xx} - 
\frac{s_{ax}^2}{s_{aa}}
\label{SSSS2}
\end{align}

Note that
the above (\ref{SSSS1}) and (\ref{SSSS2})
are the same as
(\ref{S1}).
Therefore,
Problem 2 (i.e.,
regression analysis in quantum language) is a quantum linguistic stories of the least squared method
(Problem 1).

\par
\vskip1.0cm
\par
\noindent

\subsection{Several properties (Distributions, confidence interval and hypothesis test)}
\label{subsubsec:SeveralPropperty}
\par
\noindent
\par
Since our main assertion is to mention Problem 1,
this section may be regarded as a kind of appendix.
For the detailed proofs of Lemma 1, see standard books of statistics
(e.g.,
ref.
\cite{Case}).
\par
\noindent
\par
Let
${\mathsf M}_{C_0(\Omega_{0} \times {\mathbb R}_+)}(
{\mathsf O} \equiv (X(={\mathbb R}^n), {\cal F}, F)
, S_{[(\beta_0,\beta_1, \sigma)]}
{}
)$
be the observable in Problem 3.
For each $(\beta, \sigma) \in {\mathbb R}^2 \times {\mathbb R}_+$,
we have the probability space
$(X, {\mathcal F}, P_{(\beta, \sigma)} )$,
where
$P_{(\beta, \sigma)}(\Xi ) =$
$F(\Xi)](\beta_0,\beta_1, \sigma)$
$(\forall \Xi \in {\mathcal F} )$.

Put 
$$
L^2(X)=
\{\mbox{measurable function $f:X \to {\mathbb R}$}
\;\;|\;\
 [\int_X |f(x)|^2 P_{(\beta, \sigma)}(dx)]^{1/2 } < \infty 
\}.
$$
For any $f, g \in L^2(X)$,
define $E(f)$
and
$V(f)$
such that
\begin{align}
&
E(f)= \int_X f(x) P_{(\beta, \sigma)} (dx),
\quad
V(f)=\int_X |f(x) -E(f)|^2  P_{(\beta, \sigma)} (dx).
\end{align}
\par
\noindent

\par
\noindent
\bf
Lemma 1
\sl
Consider the measurement
${\mathsf M}_{C_0(\Omega_{0} \times {\mathbb R}_+)}(
{\mathsf O} \equiv (X, {\cal F}, F)
, S_{[(\beta_0,\beta_1, \sigma)]}
{}
)$
in Problem 3.
And assume the above notations. Then, we see:
\rm
\begin{itemize}
\item[(H$_1$)]
$
\mbox{(1): }
V(\hat{\beta}_0)= \frac{\sigma^2}{n}(1+ \frac{\overline{a}^2}{s_{aa}}),
\qquad
\mbox{(2): }
V(\hat{\beta}_1)= \frac{\sigma^2}{n} \frac{1}{s_{aa}},
$
%
\item[(H$_2$)]
[Studentization].
\sl
Motivated by the (H$_1$), we see:
\begin{align}
&
T_{\beta_0}
:=
\frac{\sqrt{n}(\hat{\beta}_0-{\beta}_0)}
{\sqrt{
{\hat{\sigma}^2(1+ \overline{a}^2/ s_{aa})}}}
\sim
t_{n-2},
\qquad
T_{\beta_1}
:=
\frac{\sqrt{n}(\hat{\beta}_1-{\beta}_1)}
{\sqrt{
{\hat{\sigma}^2/ s_{aa}}}}
\sim
t_{n-2}
\end{align}
where $t_{n-2}$ is the student's distribution with $n-2$ degrees of freedom.
\end{itemize}
\rm
\rm
For the proof. see ref.
\cite{Case}.
\rm
\vskip1.0cm
Let
${\mathsf M}_{C_0(\Omega_{0}(={\mathbb R}^2) \times {\mathbb R}_+)}(
{\mathsf O} \equiv (X(={\mathbb R}^n), {\cal F}, F)
, S_{[(\beta_0,\beta_1, \sigma)]}
{}
)$
be the observable in Problem 2.
For each $k=0,1$,
define the estimator
${\widehat{E}}_k:X(={\mathbb R}^n) \to {\Theta_k}(={\mathbb R})$
and
the quantity
$\pi_k: \Omega(={\mathbb R}^2 \times {\mathbb R}_+) \to {\Theta_k}(={\mathbb R})$
as follows.
\begin{align}
&
{\widehat{E}}_0 ( x)
(=\hat{\beta}_0(x))
=
\overline{x}- \frac{s_{ax}}{s_{aa}} \overline{a},
\quad
{\widehat{E}}_1 ( x)
(=\hat{\beta}_1(x))
=
\frac{s_{ax}}{s_{aa}} ,
\quad
\pi_0
(\beta_0, \beta_1, \sigma )
=
\beta_0.
\quad
\pi_1
(\beta_0, \beta_1, \sigma )
=
\beta_1,
\\
&
\qquad \qquad \qquad
( \forall (\beta_0, \beta_1, \sigma )
\in
{\mathbb R}^2 \times {\mathbb R}_+
)
\nonumber
\end{align}


Let
$\alpha$
be a real number such that
$0 < \alpha \ll 1$,
for example,
$\alpha = 0.05$.
For any state
$ \omega =( \beta, \sigma ) 
({}\in \Omega ={\mathbb R}^2 \times {\mathbb R}_+)$,
define
the positive number
$\eta^\alpha_{\omega,  k}$
$({}> 0)$
by
(\ref{Defeta}),
that is,
\begin{align}
\eta^\alpha_{\omega, k}
&
=
\inf
\{
\eta > 0:
[F(\{ x \in X \;:\; 
d^x_{\Theta_k} ( {\widehat{E}_k}(x) , \pi_k( \omega ) )
\ge \eta
\}
)](\omega )
\le \alpha
\}
\label{eq9}
\end{align}
where, for each $\theta_k^0, \theta_k^1 (\in \Theta_k )$, the semi-distance $d_{\Theta_k}^x$ in $\Theta_k$ is defined by
\begin{align}
d^x_{\Theta_k}(\theta_k^0,\theta_k^1)
=
\begin{cases}
\frac{\sqrt{n}| \theta_0^0-\theta_0^1 |}
{\sqrt{
{\hat{\sigma}^2(1+ \overline{a}^2/ s_{aa})}}}
\quad & (\mbox{if }k=0)
\\
\\
\frac{\sqrt{n}
| \theta_1^0-\theta_1^1 |
}
{\sqrt{
{\hat{\sigma}^2/ s_{aa}}}}
\quad & (\mbox{if }k=1)
\end{cases}
\end{align}
Therefore,
we see,
by Lemma 1, that
\begin{align}
\eta^\alpha_{\omega, k}
&
=
\begin{cases}
\inf
\{
\eta > 0:
[F(\{ x \in X \;:\; 
\frac{\sqrt{n}| \hat{\beta}_0(x) - \beta_0 |}
{\sqrt{
{\hat{\sigma}^2(1+ \overline{a}^2/ s_{aa})}}}
%
%
%
%
\ge \eta
\}
)](\omega )
\le \alpha
\}
\quad & (\mbox{if }k=0)
\\
\\
\inf
\{
\eta > 0:
[F(\{ x \in X \;:\; 
\frac{\sqrt{n}|\hat{\beta}_1(x)-{\beta}_1|}
{\sqrt{
{\hat{\sigma}^2(x)/ s_{aa}}}}
\ge \eta
\}
)](\omega )
\le \alpha
\}\quad & (\mbox{if }k=1)
\end{cases}
\\
&
=
t_{n-2}(\alpha/2)
\label{Tdistri}
\end{align}

%
\vskip0.5cm
\rm
The following propositions (described in quantum language) immediately follow from (\ref{Tdistri}).
\bf
\par
\noindent
Proposition 1
\rm
[Confidence interval].
\sl
Assume that
a measured value
$x \in X$ is obtained by the measurement
${\mathsf M}_{C_0(\Omega_{0} \times {\mathbb R}_+)}(
{\mathsf O} \equiv (X, {\cal F}, F)
, S_{[(\beta_0,\beta_1, \sigma)]}
{}
)$.
Here, the state $(\beta_0,\beta_1, \sigma)$ is assumed to be unknown.
Then, we have
the
$({}1- \alpha{})$-confidence interval $I_{x,k}^{1- \alpha}$
in Theorem 2 as follows.
\begin{align}
&
I_{x,k}^{1- \alpha}
=
\{
\pi_k(\omega)
(\in
\Theta_k)
:
d^x_{\Theta_k} ({}{\widehat{E}_k}(x),
\pi_k(\omega )
)
<
\eta^{1- \alpha}_{\omega, k }
\}
\nonumber
\\
\nonumber
\\
&
=
\begin{cases}
I_{x,0}^{1- \alpha}
=
\Big\{
\beta_0 
= \pi_0(\omega)
(\in
{\Theta_0})
\;:\;
\frac{
|\hat{\beta}_0 (x) -{\beta}_0|
}{
{\sqrt{
{\frac{\hat{\sigma}^2(x)}{n}(1+ \overline{a}^2/ s_{aa})}}}
}
\le
t_{n-2}(\alpha/2)
\Big\}
\quad & (\mbox{if }k=0)
\\
\\
I_{x,1}^{1- \alpha}
=
\Big\{
\beta_1 
= \pi_1(\omega)
(\in
{\Theta_1})
:
\frac{
|\hat{\beta}_1 (x) -{\beta}_1|
}{
{\sqrt{
{\frac{\hat{\sigma}^2(x)}{n}(1/ s_{aa})}}}
%
}
\le
t_{n-2}(\alpha/2) 
\Big\}
\quad & (\mbox{if }k=1)
%
%
%
\end{cases}
\label{eq12} 
\end{align}

\par
\vskip1.0cm
\par
\noindent
\bf
Proposition 2
\rm
[Hypothesis test].
\sl
Consider the measurement
${\mathsf M}_{C_0(\Omega_{0} \times {\mathbb R}_+)}(
{\mathsf O} \equiv (X, {\cal F}, F)
, S_{[(\beta_0,\beta_1, \sigma)]}
{}
)$.
Here, the state $(\beta_0,\beta_1, \sigma)$ is assumed to be unknown.
Then, according to Theorem 2, we say:
\begin{itemize}
\rm
\item[(I$_1$)]
\sl
Assume the null hypothesis
$H_{N} = { \{ \beta_0 \}}
(\subseteq \Theta_0={\mathbb R})$.
Then, the rejection region
is as follows:
\begin{align}
{\widehat R}_{{H_N}}^{\alpha; X}
&
=
{\widehat{E}_0}^{-1}(
{\widehat R}_{{H_N}}^{\alpha; {\Theta_0}})
=
\bigcap_{\omega \in  \Omega \mbox{ \footnotesize such that }
\pi_0(\omega) \in {H_N}}
\{
x
(\in
X)
:
d^x_{\Theta_0} ({}{\widehat{E}_0}(x),
\pi_0(\omega )
)
\ge
\eta^\alpha_{\omega }
\}
\nonumber
\\
&
=
\Big\{
x \in X
\;:\;
\frac{
|\hat{\beta}_0 (x) -{\beta}_0|
}{
{\sqrt{
{\frac{\hat{\sigma}^2(x)}{n}(1+ \overline{a}^2/ s_{aa})}}}
}
\ge
t_{n-2}(\alpha/2) 
\Big\}
\end{align}
\rm
\item[(I$_2$)]
\sl
Assume the null hypothesis
$H_N = { \{ \beta_1 \}}
(\subseteq \Theta_1={\mathbb R})$.
Then, the rejection region
is as follows:
\begin{align}
{\widehat R}_{{H_N}}^{\alpha; X}
&
=
{\widehat{E}_1}^{-1}(
{\widehat R}_{{H_N}}^{\alpha; {\Theta_1}})
=
\bigcap_{\omega \in  \Omega \mbox{ \footnotesize such that }
\pi_1(\omega) \in {H_N}}
\{
x
(\in
X)
:
d^x_{\Theta_1} ({}{\widehat{E}_1}(x),
\pi_1(\omega )
)
\ge
\eta^\alpha_{\omega }
\}
\nonumber
\\
&
=
\Big\{
x \in X
\;:\;
\frac{
|\hat{\beta}_1 (x) -{\beta}_1|
}{
{\sqrt{
{\frac{\hat{\sigma}^2(x)}{n}(1/ s_{aa})}}}
%
}
\ge
t_{n-2}(\alpha/2) 
\Big\}
\end{align}
\end{itemize}

\par
\vskip1.0cm
\par
\subsection{The quantum linguistic formulation of generalized linear model}
\label{subsec:23GLM}
\rm

As the generalization of Section \ref{subsec:21GLM},
we shall discuss
the generalized linear model in quantum language as follows:

\par
\noindent
\par
Put
$T=\{ 0,1,2, \cdots, i , \cdots, n \}$,
which is the same as the tree (\ref{Tree}),
that is,
\begin{align}
\tau(i)=0
\qquad
(\forall i =1,2, \cdots, n)
\end{align}
For each $i \in T$,
define a locally compact space $\Omega_i$ such that
\begin{align}
&
\Omega_{0}={\mathbb R}^{m+1}
=
\Big\{
\beta
=\begin{bmatrix}
\beta_0 \\
\beta_1 \\
\vdots
\\
\beta_m
\end{bmatrix}
\;:\;
\beta_0, \beta_1, \cdots, \beta_m \in {\mathbb R}
\Big\}
\quad
\\
&
\Omega_{i}={\mathbb R}
=
\Big\{
\mu_i
\;:\;
\mu_i \in {\mathbb R}
\Big\}
\quad(i=1,2, \cdots, n )
\end{align}
Assume that
\begin{align}
a_{ij} \in {\mathbb R} 
\qquad
(i=1,2, \cdots, n, \;\;j=1,2, \cdots, m, (m+1 \le n) )
\label{Explanatory2},
\end{align}
which are called {\it explanatory variables}
in the conventional statistics.
Consider the deterministic causal map
$\psi_{a_{i \tiny{\bullet}}}: \Omega_0(={\mathbb R}^{m+1})
\to
\Omega_{i} (={\mathbb R})$
such that
\begin{align}
&
\Omega_0={\mathbb R}^{m+1} \ni \beta =(\beta_0, \beta_1, \cdots, \beta_m ) \mapsto
\psi_{a_{i \tiny{\bullet}}} ( \beta_0, \beta_1, \cdots, \beta_m )
=
\beta_0 + \sum_{j=1}^m \beta_j a_{ij}= \mu_i \in \Omega_i ={\mathbb R}
\\
&
\qquad \qquad \qquad \qquad
\qquad \qquad \qquad \qquad
(i=1,2, \cdots, n)
\nonumber
\end{align}
Summing up, we see
\begin{align}
\beta=
\begin{bmatrix}
\beta_0 \\
\beta_1 
\\
\beta_2 \\
\vdots
\\
\beta_m
\end{bmatrix}
\mapsto
\begin{bmatrix}
\psi_{a_{1 \tiny{\bullet}}} ( \beta_0, \beta_1, \cdots, \beta_m) \\
\psi_{a_{2 \tiny{\bullet}}} ( \beta_0, \beta_1, \cdots, \beta_m) 
\\
\psi_{a_{3 \tiny{\bullet}}} ( \beta_0, \beta_1, \cdots, \beta_m) \\
\vdots
\\
\psi_{a_{n \tiny{\bullet}}} ( \beta_0, \beta_1, \cdots, \beta_m)
\end{bmatrix}
=
\begin{bmatrix}
 1 & a_{1{1}} & a_{12} & \cdots & a_{1m} 
 \\
 1 & a_{2{1}} & a_{22} & \cdots & a_{2m} 
 \\
  1 & a_{3{1}} & a_{32} & \cdots & a_{3m} 
  \\
  1 & a_{4{1}} & a_{42} & \cdots & a_{4m} 
  \\
\vdots & \vdots & \vdots & \vdots &  \vdots 
\\
 1 & a_{n{1}} & a_{n2} & \cdots & a_{nm}
\end{bmatrix}
\cdot
\begin{bmatrix}
\beta_0 \\
\beta_1 
\\
\beta_2 \\
\vdots
\\
\beta_m
\end{bmatrix}
\label{GeneralLinear444}
\end{align}
which is equivalent to the deterministic Markov operator
$\Psi_{a_{i \tiny{\bullet} }}:
C_0(\Omega_{i}) \to C_0(\Omega_0)$
such that
\begin{align}
[{\Psi_{a_{i \tiny{\bullet} }}}(f_i)](\omega_0)
=
f_i( \psi_{a_{i \tiny{\bullet} }} (\omega_0))
\quad
(\forall f_i \in C_0(\Omega_{i}),
\;\;
\forall \omega_0 \in \Omega_0,
\forall i \in 1,2, \cdots, n)
\end{align}
Thus, under the identification: $a_{ij} \Leftrightarrow \Psi_{a_{i \tiny{\bullet}}}$,
the term "explanatory variable" means a kind of causality.
\par
\noindent
\par
\noindent
\begin{figure}[htbp]
{
\setlength{\unitlength}{0.7mm}
\begin{picture}(80,55)(0,20)
\put(40,0){
\put(83,60){\makebox(10,10)[r]{${C_0 (\Omega_1(\equiv{\mathbb R}))}$}}
\put(83,45){\makebox(10,10)[r]{${C_0 (\Omega_2(\equiv{\mathbb R}))}$}}
\put(83,22){\makebox(10,10)[r]{${C_0 (\Omega_n(\equiv{\mathbb R}))}$}}
\put(29,41){\makebox(10,10)[r]{${C_0 (\Omega_0(\equiv{\mathbb R}^{m+1}))}$}}
\put(60,65){\vector(-3,-2){15}}
\put(60,50){\vector(-3,-1){13}}
\put(60,27){\vector(-3,2){15}}
%
\put(47,67){$\Psi_{a_{1\tiny{\bullet}}}$}
\put(50,51){$\Psi_{a_{2 \tiny{\bullet} }}$}
\put(50,43){$\cdots \cdots$}
\put(50,38){$\cdots \cdots$}
\put(47,25){$\Psi_{a_{n \tiny{\bullet} } }$}
}
\end{picture}
}
\caption{
\label{Fig3}
Parallel structure(Causal relation $\Psi_{a_{i \tiny{\bullet}}}$)
}
\end{figure}
\noindent
Therefore, we have
the observable
${\mathsf O}_{0}^{a_{i \tiny{\bullet} }} {{\equiv}} ({\mathbb R}, {\cal B}_{{\mathbb R}}, \Psi_{a_{i \tiny{\bullet} }}G_{\sigma_{}})$
in
$C_0 (\Omega_{0}(\equiv
{\mathbb R}^{m+1}))$
such that
\begin{align}
&
[\Psi_{a_{i \tiny{\bullet} }}(G_{\sigma_{}}(\Xi))] (\beta ) = 
[(G_{\sigma_{}}(\Xi))] (\psi_{a_{i \tiny{\bullet} }}(\beta ))
=
\frac{1}{(\sqrt{2 \pi \sigma_{}^2})}
\underset{\Xi}{\int} \exp
\Big[{- \frac{   (x - (\beta_0 + \sum_{j=1}^m a_{i{j}} \beta_j ))^2}{2 \sigma_{}^2}}
\Big] dx
\\
&
\qquad
\qquad
(\forall \Xi \in {\cal B}_{{\mathbb R}}, \forall 
\beta
=(\beta_0, \beta_1, \cdots, \beta_m )\in \Omega_{0} (\equiv {\mathbb R}^{m+1} 
))
\nonumber
\end{align}
\par
\noindent
\par
\noindent
Hence, we have the simultaneous observable
$\bigtimes_{i=1}^n{\mathsf O}_{0}^{a_{i \tiny{\bullet} }} {{\equiv}} ({\mathbb R}^n, {\cal B}_{{\mathbb R}^n}, \bigtimes_{i=1}^n \Psi_{a_{i \tiny{\bullet} }}G_{\sigma_{}})$
in
$C_0 (\Omega_{0}(\equiv
{\mathbb R}^{m+1}))$
%
such that
\begin{align}
&
[(\bigtimes_{i=1}^n \Psi_{a_{i \tiny{\bullet} }}G_{\sigma_{}})
(\bigtimes_{i=1}^n \Xi_i)](\beta)
=
\bigtimes_{i=1}^n \Big(
[\Psi_{a_{i \tiny{\bullet} }}G_{\sigma_{}})
(\Xi_i)](\beta)\Big)
\nonumber
\\
=
&
\frac{1}{(\sqrt{2 \pi \sigma_{}^2})^n}
\underset{\bigtimes_{i=1}^n  \Xi_i}{\int \cdots \int} \exp
\Big[{- \frac{ \sum_{i=1}^n  (x_i - (\beta_0 + \sum_{j=1}^ma_{i{j}} \beta_j ))^2}{2 \sigma_{}^2}}
\Big] dx_1 \cdots dx_n
\label{Fundamental44}
\\
&
\qquad
\qquad
(\forall \bigtimes_{i=1}^n \Xi_i \in {\cal B}_{{\mathbb R}^n}, \forall 
\beta
=(\beta_0, \beta_1, \cdots, \beta_m )
\in \Omega_{0} (\equiv {\mathbb R}^{m+1} ))
\nonumber
\end{align}
Assuming that $\sigma$ is variable, we have the observable
${\mathsf O}=
\Big({\mathbb R}^n(=X) , {\mathcal B}_{{\mathbb R}^n}(={\mathcal F}),
F \Big)$ in $C_0 ( \Omega_0 \times {\mathbb R}_+ )$
such that
\begin{align}
[F(\bigtimes_{i=1}^n \Xi_i )](\beta, \sigma_{})
=
[(\bigtimes_{i=1}^n \Psi_{a_{i \tiny{\bullet} }}G_{\sigma_{}})
(\bigtimes_{i=1}^n \Xi_i)](\beta)
\quad
(\forall \bigtimes_{i=1}^n \Xi_i \in {\cal B}_{{\mathbb R}^n},
\forall
(\beta , \sigma_{} ) \in {\mathbb R}^{m+1} (\equiv \Omega_0) \times {\mathbb R}_+ )
\end{align}

\par


\par
\noindent
\bf
Problem 3
\rm
[The generalized linear model]
\sl
Assume that
a measured value 
$x=\begin{bmatrix}
x_1
\\
x_2
\\
\vdots
\\
x_n
\end{bmatrix}
\in X={\mathbb R}^n $
is obtained by the measurement
${\mathsf M}_{C_0(\Omega_{0} \times {\mathbb R}_+)}(
{\mathsf O} \equiv (X, {\cal F}, F)
, S_{[(\beta_0,\beta_1,\cdots, \beta_m , \sigma)]}
{}
)$.
We do not know the state
$
(\beta_0, \beta_1,\cdots, \beta_m, \sigma_{}^2
)$.
%
%
%
Then, from the measured value
$
x=(
x_1, x_2, \ldots, x_n )
\in {\mathbb R}^{n}$,
infer the $\beta_0, \beta_1, \cdots, \beta_m, \sigma_{}$!
That is,
represent the
$(\beta_0, \beta_1, \cdots, \beta_m, \sigma_{})$
by $(\hat{\beta}_0(x), \hat{\beta}_1(x),\cdots, \beta_m(x), \hat{\sigma}_{}(x))$
(i.e.,
the functions of $x$).
\rm
\vskip0.5cm
\par
\noindent
\bf
Answer.
\rm
The answer is easy, since it is a slight generalization of Problem 2. Also, it suffices to follow
ref. \cite{Case}.
However, note that the purpose of this paper is 
to describe Problem 3
(i.e,
the quantum linguistic formulation of
the generalized linear model)
and not to give the answer to Problem 3.
\vskip0.5cm
\par

\rm
\section{Conclusions
}
\par
\noindent
\par
\rm
\par
\noindent
\par
Quantum language is clearly defined by the (B), that is,
\begin{itemize}
\item[(B$_1$)]$\underset{\mbox{(=MT(measurement theory))}}{\fbox{Quantum language}}
=
\underset{\mbox{(measurement)}}{\fbox{Axiom 1}}
+
\underset{\mbox{(causality)}}{\fbox{Axiom 2}}
+
\underset{\mbox{(how to use Axioms)}}{\fbox{linguistic interpretation}}
\label{eq1}
$
\end{itemize}
Therefore, we do not start from "random variable" but "measurement".
Our purpose of this paper was to understand the regression analysis and
the generalized linear model in quantum language.
In fact, we showed
\begin{itemize}
\item[(J)]
the term "explanatory variable in (\ref{Explanatory1})
and
(\ref{Explanatory2})" is characterized a kind of causality
(cf. Figures \ref{Fig2} and \ref{Fig3}).
And the term "response variable" means the measured value.
\end{itemize}

We believe that quantum language has a great power of
description,
and therefore, even statistics can be described by quantum language.
We hope that our assertions will be examined from various points of view.

%
%



\rm
\par
\renewcommand{\refname}{
\large 
References}
{
\small

\normalsize
}
\end{document}